\documentclass[brochure,12pt]{bourbaki}
\usepackage[matrix,arrow]{xy}
\usepackage{amssymb,amsfonts,amsmath,footnote}
\usepackage[latin1]{inputenc}
\usepackage[T1]{fontenc}
\usepackage{lmodern}
\usepackage[french]{babel}
\date{juin 2015}
\bbkannee{67\`eme ann\'ee, 2014-2015}
\bbknumero{1102}
\title{Construction de repr\'esentations galoisiennes\\ de torsion}
\subtitle{d'apr\`es Peter Scholze}
\author{Sophie MOREL}
\address{Department of Mathematics\\ 
Princeton University\\ 
Fine Hall\\
 Washington Road \\
Princeton, NJ 08544, U.S.A.} 
\email{smorel@math.princeton.edu}

\usepackage{hyperref}

\usepackage{mathrsfs}

\newcommand{\C}{\mathbb{C}}
\newcommand{\Nat}{\mathbb{N}}
\newcommand{\R}{\mathbb{R}}
\newcommand{\Q}{\mathbb{Q}}
\newcommand{\Z}{\mathbb{Z}}
\newcommand{\Fi}{\mathbb{F}}

\newcommand{\Proj}{\mathbb{P}}
\newcommand{\Ade}{\mathbb{A}}
\newcommand{\Af}{{\mathbb{A}_f}}

\renewcommand{\H}{{H}}

\newcommand{\GL}{{\bf{GL}}}

\newcommand{\SL}{\bf{SL}}
\newcommand{\GSp}{{\bf{GSp}}}
\newcommand{\Sp}{{\bf {Sp}}}

\newcommand{\SO}{\bf{SO}}

\newcommand{\Pa}{{\bf P}}
\newcommand{\QP}{{\bf Q}}

\newcommand{\N}{{\bf N}}

\newcommand{\Se}{{\bf S}}
\newcommand{\Ar}{\mathrm{A}}

\newcommand{\G}{{\bf{G}}}
\newcommand{\GH}{{\bf{H}}}

\newcommand{\Of}{\mathcal{O}}
\newcommand{\Hf}{\mathcal{H}}

\newcommand\Abf{\mathcal{A}}

\newcommand{\et}{\mathrm{\acute{e}t}} %% plus accent
\newcommand{\proet}{\mathrm{pro\acute{e}t}} %% plus accent

\newcommand{\fl}{\longrightarrow}
\newcommand{\Flr}{\mathrm{Fl}}
\newcommand{\Fl}{{\mathscr{F}\!\ell}}

\newcommand{\fle}{\longmapsto}

\DeclareMathOperator{\Gal}{Gal} %% for Galois group
\newcommand\ggoth{\mathfrak{g}}
\newcommand\Ha{\mathrm{Ha}}
\newcommand{\Hecke}{\mathcal{H}}
\newcommand\Ical{\mathcal{I}}

\newcommand{\K}{\mathrm{K}}
\DeclareMathOperator{\Ker}{Ker}
\DeclareMathOperator{\Lie}{Lie}
 
\newcommand{\Lf}{\mathcal{L}} %% syst�me local
\newcommand{\M}{\mathcal{M}} %% cat�gorie de faisceaux pervers

\mathchardef\mhyphen="2D

\newcommand\oQ{\overline{\Q}}

\newcommand\quash[1]{}
\newcommand{\sous}{\setminus}

\DeclareMathOperator{\Spa}{Spa}

\newcommand{\ungras}{1\!\!\mkern -1mu1}
\newcommand{\Vcal}{\mathcal{V}}
\newcommand{\X}{\mathcal{X}}

\newcommand{\Y}{\mathcal{Y}}

\begin{document}
\maketitle

\noindent{\bf INTRODUCTION}

\bigskip

Le but de cet expos\'e est de pr\'esenter les r\'esultats de Scholze sur
la construction des repr\'esentations galoisiennes de torsion associ\'ees
aux caract\`eres de l'alg\`ebre de Hecke apparaissant dans la cohomologie de
torsion des espaces localement sym\'etriques associ\'es au groupe $\GL_n$.
La phrase pr\'ec\'edente
est expliqu\'ee plus en d\'etail dans la section 1.
Toutes les erreurs et inexactitudes dans ce texte sont 
bien entendu dues \`a l'auteur et
non \`a Scholze.

Je remercie Ana Caraiani pour d'utiles remarques sur une version pr\'ec\'edente
de cet expos\'e.

\section{Quelques conjectures du programme de Langlands}
\subsection{Repr{\'e}sentations automorphes}

La r{\'e}f{\'e}rence standard pour la d{\'e}finition des formes et repr{\'e}sentations
automorphes est le texte de Borel et Jacquet dans Corvallis (\cite{BJ}), voir
aussi le livre de Moeglin et Waldspurger (\cite{MW}).

Soit $n$ un entier strictement positif. On note $\Ade=\R\times\Af$ l'anneau des
ad{\`e}les de $\Q$, o{\`u}
\[\Af=\Q\otimes_\Z\widehat{\Z}=\{(x_p)\in\prod_{p\ {\rm premier}}\Q_p|x_p\mbox{ pour
presque tout }p\}\]
est l'anneau des ad\`eles finies
(o{\`u} {\og pour presque tout $p$\fg} signifie {\og pour tout $p$ sauf un nombre fini\fg}).
On fixe une mesure de Haar sur le groupe topologique
\[\GL_n(\Ade)=\{(g_\infty,(g_p))\in\GL_n(\R)\times\prod_p\GL_n(\Q_p)|g_p\in\GL_n(
\Z_p)\mbox{ pour presque tout }p\},\]
et on note $L^2_{\GL_n}=L^2(\GL_n(\Q)\Ar\sous\GL_n(\Ade),\C)$,
o{\`u} $\GL_n(\Q)$ est plong{\'e}
diagonalement dans $\GL_n(\Ade)$ et $\Ar=\R_{>0}$ est la composante connexe
de $1$ dans le centre de $\GL_n(\R)$.\footnote{On pourrait fixer un caract{\`e}re
unitaire quelconque $\xi$ de $\Ar$ et consid{\'e}rer l'espace des fonctions
$f:\GL_n(\Q)\sous\GL_n(\Ade)\fl\C$ telles que $f(zg)=\xi(z)f(g)$ pour tous
$z\in\Ar$ et $g\in\GL_n(\Ade)$ et qui sont $L^2$ modulo $\Ar$. Dans cet expos{\'e},
on prendra $\xi=1$, mais c'est uniquement pour all{\'e}ger les notations.}

Le groupe $\GL_n(\Ade)$ agit sur $L^2_{\GL_n}$ par translation {\`a} droite sur
l'argument de la fonction. Une \emph{repr{\'e}sentation automorphe discr{\`e}te de
$\GL_n(\Ade)$} est une repr{\'e}sentation irr{\'e}ductible de $\GL_n(\Ade)$ qui
appara{\^i}t comme facteur direct de la repr{\'e}sentation $L^2_{\GL_n}$. En fait, on a
\[L^2_{\GL_n}=L^2_{\GL_n,{\rm disc}}\oplus L^2_{\GL_n,{\rm cont}}\]
en tant que repr{\'e}sentation de $\GL_n(\Ade)$,
o{\`u} $L^2_{\GL_n,{\rm cont}}$ n'a pas de facteur direct irr{\'e}ductible et
\[L^2_{\GL_n,disc}=\bigoplus \pi^{m(\pi)}\]
(somme directe compl{\'e}t{\'e}e),
o{\`u} la somme est sur les repr{\'e}sentations automorphes discr{\`e}tes et les $m(\pi)$
sont des entiers strictement positifs.

Soit $f\in L^2_{\GL_n}$ une fonction born{\'e}e. On dit que $f$ est \emph{cuspidale}
si, pour tout sous-groupe parabolique propre $\Pa$ de $\GL_n$, si on note
$\N_P$ le radical unipotent de $\Pa$
et $dn$ une mesure de Haar sur $\N_P(\Ade)$,
alors, pour tout $g\in\GL_n(\Ade)$,
\[\int_{\N_P(\Q)\sous\N_P(\Ade)} f(ng)dn=0.\]
L'espace $L^2_{\GL_n,{\rm cusp}}$ des fonctions born{\'e}es cuspidales est un sous-espace
de $L^2_{\GL_n}$ ferm{\'e} et stable par l'action de $\GL_n(\Ade)$,
qui est contenu dans $L^2_{\GL_n,{\rm disc}}$, c'est-{\`a}-dire
somme directe compl{\'e}t{\'e}e de repr{\'e}sentations irr{\'e}ductibles de $\GL_n(\Ade)$.
Les repr{\'e}sentations irr{\'e}ductibles de $\GL_n(\Ade)$
qui apparaissent dans $L^2_{\GL_n,{\rm cusp}}$
sont dites \emph{automorphes cuspidales}.

De plus, si $\pi$ est une repr{\'e}sentation automorphe discr{\`e}te de $\GL_n(\Ade)$,
on a
\[\pi=\pi_\infty\otimes\bigotimes'_{p\ premier}\pi_p,\]
o{\`u} $\pi_\infty$ (resp. $\pi_p$) est une repr{\'e}sentation irr{\'e}ductible de
$\GL_n(\R)$ (resp. $\GL_n(\Q_p)$) et, pour presque tout $p$, la repr{\'e}sentation
$\pi_p$ est \emph{non ramifi{\'e}e}, c'est-{\`a}-dire que $\pi_p^{\GL_n(\Z_p)}\not=0$
(cet espace est alors de dimension $1$). Voir l'article de Flath
\cite{F} pour la
d{\'e}finition du produit tensoriel restreint $\bigotimes'$ et pour des r{\'e}f{\'e}rences.
La classification de Langlands associe {\`a} la repr{\'e}sentation irr{\'e}ductible
admissible $\pi_\infty$ de $\GL_n(\R)$ une repr{\'e}sentation du groupe de
Weil $W_\R$ de $\R$ dans $\GL_n(\C)$. En restreignant cette repr{\'e}sentation
au sous-groupe $\C^\times$ de $W_\R$, on obtient un morphisme $r:\C^\times\fl
\GL_n(\C)$. On dit que la repr{\'e}sentation automorphe $\pi$ est \emph{alg{\'e}brique}
si $r$ est un morphisme de groupes alg{\'e}briques sur $\C$.
Cette d{\'e}finition est due {\`a} Clozel (d{\'e}finition 1.8 de \cite{C1}), et peut aussi
se formuler comme une condition d'int{\'e}gralit{\'e} sur le caract{\`e}re infinit{\'e}simal de $\pi_\infty$, c'est-{\`a}-dire
le caract{\`e}re par lequel
le centre de l'alg{\`e}bre universellement enveloppante de $\Lie(\GL_n(\R))\otimes_
\R\C$ agit sur $\pi_\infty$.
\footnote{En fait, la d{\'e}finition donn{\'e}e ci-dessus n'est pas tout {\`a} fait
celle de Clozel, car on a supprim{\'e} la torsion par $|.|^{(n-1)/2}$. La notion
que nous avons d{\'e}finie est celle de repr{\'e}sentation \emph{L-alg{\'e}brique}
au sens de Buzzard et Gee (cf. \cite{BG}), qui semble plus adapt{\'e}e au cas
d'un groupe r{\'e}ductif g{\'e}n{\'e}ral, et qui est celle qu'utilise Scholze.}

\subsection{Conjecture de r{\'e}ciprocit{\'e} de Langlands et Clozel}

\begin{conj}\label{conj:LC}
Soit $\pi=\pi_\infty\otimes\bigotimes'_p \pi_p$ une repr{\'e}sentation
automorphe alg{\'e}brique cuspidale de $\GL_n(\Ade)$, et soit $\ell$ un nombre
premier.
Alors il existe une repr{\'e}sentation continue semi-simple\footnote{
En fait, la repr\'esentation $\rho_\pi$ devrait m\^eme \^etre irr\'eductible,
mais on ne sait le prouver que dans quelques cas particuliers, voir par
exemple \cite{BR} et le th\'eor\`eme D de \cite{BLGGT}.}
$\rho_\pi:\Gal(\overline{\Q}/\Q)\fl\GL_n(\oQ_\ell)$ telle que, pour tout
nombre premier $p\not=\ell$ tel que $\pi_p$ soit non ramifi{\'e}e, $\pi_p$
et $\rho_{\pi|\Gal(\overline{\Q}_p/\Q_p)}$ se correspondent par l'isomorphisme
de Satake.

\end{conj}

Donnons quelques explications sur l'{\'e}nonc{\'e}. Une r{\'e}f{\'e}rence pour l'isomorphisme
de Satake est le chapitre IV de l'article \cite{C} de Cartier, voir aussi
l'article introductif \cite{G} de Gross. Soit $p$ un nombre
premier. Rappelons que l'on dit que $\pi_p$ est non ramifi{\'e}e (ou que
$\pi$ est non ramifi{\'e}e en $p$) si $\pi_p^{\GL_n(\Z_p)}\not=0$. Cet espace
d'invariants est alors n{\'e}cessairement de dimension $1$, et d{\'e}finit donc
un caract{\`e}re de l'alg{\`e}bre de Hecke non ramifi{\'e}e $\Hecke_p$ des fonctions
$f:\GL_n(\Q_p)\fl\C$ {\`a} support compact et invariantes {\`a} gauche et {\`a} droite
par $\GL_n(\Z_p)$ (le produit est le produit de convolution, d{\'e}fini en
utilisant la mesure de Haar sur $\GL_n(\Q_p)$ telle que $\GL_n(\Z_p)$ soit
de volume $1$). Il r{\'e}sulte de l'isomorphisme de Satake (sous la forme, par
exemple, du corollaire 4.2 de \cite{C}) que l'ensemble des caract{\`e}res
de $\Hecke_p$ est en bijection avec l'ensemble des classes de conjugaison
semi-simples de $\GL_n(\C)$.

Revenant {\`a} la conjecture, on note $a_{\pi_p}$ la classe de conjugaison
correspondant {\`a} $\pi_p$. Dire que $\rho_{\pi|\Gal(\overline{\Q}_p/\Q_p)}$ correspond
{\`a} $\pi_p$ par l'isomorphisme de Satake signifie d'abord que
$\rho_{\pi|\Gal(\overline{\Q}_p/\Q_p)}$ est non ramifi{\'e}e, c'est-{\`a}-dire se
factorise par le quotient $\Gal(\overline{\Fi}_p/\Fi_p)$ de
$\Gal(\overline{\Q}_p/\Q_p)$, et ensuite que l'image par $\rho_\pi$ du morphisme
de Frobenius g{\'e}om{\'e}trique (i.e. le g{\'e}n{\'e}rateur
$x\fle x^{1/p}$ de $\Gal(\overline{\Fi}_p/\Fi_p)$) est dans la classe de
conjugaison $a_{\pi_p}$.\footnote{Le plongement $\Gal(\overline{\Q}_p/\Q_p)
\fl\Gal(\overline{\Q}/\Q)$ n'est pas canonique, puisqu'il d{\'e}pend du choix
d'un plongement $\overline{\Q}\fl\overline{\Q}_p$. Cependant, la d{\'e}finition
ci-dessus ne d{\'e}pend pas de ce choix.}

D'apr{\`e}s le th{\'e}or{\`e}me de densit{\'e} de {\v C}eboratev, $\rho_\pi$ est uniquement
d{\'e}termin{\'e}e par $\pi$. D'apr{\`e}s le th{\'e}or{\`e}me de multiplicit{\'e} un fort
de Piatetski-Shapiro et Jacquet-Shalika (cf. \cite{PS}), $\pi$
est uniquement d{\'e}termin{\'e}e par $\rho_\pi$.

\begin{rema} En combinant les conjectures de Langlands, Clozel et Fontaine-Mazur, on obtient en fait une bijection conjecturale entre les classes d'isomorphisme
de repr{\'e}sentations automorphes cuspidales alg{\'e}briques de $\GL_n(\Ade)$ et les
classes d'isomorphisme de repr{\'e}sentations continues irr{\'e}ductibles
$\Gal(\oQ/\Q)\fl\GL_n(\oQ_\ell)$ qui sont g{\'e}om{\'e}triques (c'est-{\`a}-dire
presque partout non ramifi{\'e}es et de Rham en $\ell$, voir \cite{FM}).

\end{rema}

\begin{rema} En fait, Langlands conjecture qu'il existe un groupe
pro-alg{\'e}brique $\Lf_\Q$ sur $\C$ (le groupe de Langlands de $\Q)$ tel que,
pour tout entier $n$, les repr{\'e}sentations
alg\'ebriques irr\'eductibles de dimension $n$ de $\Lf_\Q$
classifient les repr{\'e}sentations
automorphes cuspidales de $\GL_n(\Ade)$.
(En d'autres termes, les repr{\'e}sentations automorphes
cuspidales de $\GL_n(\Ade)$ sont les objets simples de dimension $n$
d'une cat{\'e}gorie tannakienne, dont le groupe tannakien est $\Lf_\Q$.)
Le groupe de Galois motivique $\M_\Q$
de $\Q$ serait alors le quotient de $\Lf_\Q$
correspondant {\`a} la sous-cat{\'e}gorie des repr{\'e}sentations alg{\'e}briques. Si
$\pi$ est une repr{\'e}sentation automorphe cuspidale alg{\'e}brique de $\GL_n(\Ade)$
et $\varphi:\Lf_\Q\fl\M_\Q
\fl\GL_{n,\C}$ est la repr{\'e}sentation alg{\'e}brique correspondante,
on devrait obtenir $\rho_\pi$ en {\'e}valuant $\varphi$ sur les $\oQ_\ell$-points
(on choisit un isomophisme $\oQ_\ell\simeq\C$ pour faire ceci), puis en
restreignant au sous-groupe $\Gal(\overline{\Q}/\Q)$ de $\M_\Q(\oQ_\ell)$
(le plongement {\'e}tant donn{\'e} par la r{\'e}alisation $\ell$-adique).

\end{rema}

\subsection{Repr{\'e}sentations automorphes et cohomologie des espaces localement
sym{\'e}triques}

Toutes les preuves de cas particuliers de la conjecture \ref{conj:LC}\footnote{
Au moins celles connues de l'auteur.} passent pas l'{\'e}tude de la cohomologie
des espaces localement sym{\'e}triques. Si $K$ est un sous-groupe compact ouvert
de $\GL_n(\Af)$ (par exemple un sous-groupe d'indice fini de $\GL_n(\widehat
{\Z})$), on pose
\[X_K=\GL_n(\Q)\sous\GL_n(\Ade)/K \Ar K_\infty,\]
o{\`u} $K_\infty=\SO(n)\subset\GL_n(\R)$. C'est une vari{\'e}t{\'e} analytique r{\'e}elle si
$K$ est assez petit (sinon, c'est un orbifold). Dans la suite, on supposera
toujours $K$ assez petit.

Si $K$ et $K'$ sont deux sous-groupes ouverts compacts de $\GL_n(\Af)$ et
$g\in\GL_n(\Af)$ est tel que $g^{-1}K'g\subset K$, on a un morphisme
analytique fini $c_g:X_{K'}\fl X_K$ qui envoie la classe de $h\in\GL_n(\Ade)$ sur
celle de $hg$. Donc, si on prend $K'=K\cap gK g^{-1}$, on obtient une
correspondance $(c_g,c_1):X_{K'}\fl X_K\times X_K$ (appel{\'e}e
\emph{correspondance de Hecke}), qui induit un morphisme
$u_g:H^*_c(X_K)\fl H^*_c(X_K)$, o{\`u} $H^*(X_K)$ est la cohomologie de Betti {\`a}
supports compacts et {\`a} coefficients dans $\C$ de $X_K$ (voir \cite{P} 1.2, 1.3).

Soit $\Hf$ l'alg{\`e}bre de Hecke globale, c'est-{\`a}-dire l'alg{\`e}bre des fonctions
localement constantes {\`a} support compact de $\GL_n(\Af)$ dans $\C$, munie du
produit de convolution (pour une mesure de Haar fix{\'e}e sur $\GL_n(\Af)$). C'est
une alg{\`e}bre associative non unitaire. Rappelons qu'une repr{\'e}sentation
$\pi:\GL_n(\Af)\fl\GL(V)$ (o{\`u} $V$ est un $\C$-espace vectoriel) est appel{\'e}e
\emph{lisse} si $\pi$ est continue pour la topologie discr{\`e}te sur $\GL(V)$,
et \emph{admissible} si, pour tout sous-groupe compact ouvert $K$ de
$\GL_n(\Af)$, $V^K$ est de dimension finie. Par exemple, si $\pi=\pi_\infty\otimes\pi_f$ est une repr{\'e}sentation automorphe discr{\`e}te de $\GL_n(\Ade)=
\GL_n(\R)\times\GL_n(\Af)$, alors $\pi_f$ est lisse admissible.

La cat{\'e}gorie des repr{\'e}sentations lisses de $\GL_n(\Af)$ est naturellement
{\'e}quivalente {\`a} celle des repr{\'e}sentations $V$ de $\Hf$ telles que, pour tout
$v\in V$, on ait $\ungras_K.v=v$ pour $K\subset\GL_n(\Af)$ assez petit. En
particulier, toute repr{\'e}sentation automorphe discr{\`e}te de $\GL_n(\Ade)$
d\'efinit une repr{\'e}sentation de $\Hf$.

D'autre part, on d{\'e}finit une action de $\Hf$ sur $H^*_c:=\varinjlim_K H^*_c(X_K)$
en convenant que, pour tout sous-groupe compact ouvert $K$ de $\GL_n(\Af)$
et tout $g\in\GL_n(\Af)$, la fonction caract{\'e}ristique de $Kg^{-1}$ agit
sur $H^*_c(X_K)$ par l'op{\'e}rateur $u_g$.

Enfin, on dit qu'une repr{\'e}sentation automorphe discr{\`e}te $\pi=\pi_\infty\otimes
\pi_f$ est \emph{cohomologique} s'il existe une repr{\'e}sentation alg{\'e}brique
$W$ de $\GL_n(\R)$ telle que la $(\ggoth,\Ar\K_\infty)$-cohomologie de
$\pi_\infty\otimes W$ soit non nulle (o{\`u} $\ggoth=\Lie(\GL_n(\R))$).
Cela implique
que $\pi$~est alg{\'e}brique et que le caract{\`e}re infinitis{\'e}mal de $\pi_\infty$
v{\'e}rifie de plus une condition de r{\'e}gularit{\'e}.\footnote{Les
$\pi_\infty$ possibles ont
{\'e}t{\'e} classifi{\'e}es par Vogan et Zuckerman dans \cite{VZ}.}

Il r{\'e}sulte alors de la conjecture de Borel, prouv{\'e}e par Franke (th{\'e}or{\`e}me
18 de \cite{Fr}), que l'on a le th{\'e}or{\`e}me suivant :

\begin{theo}\label{th:Franke}
Les sous-quotients irr{\'e}ductibles de la repr{\'e}sentation de
$\Hf$ sur $H^*_c$ viennent tous de repr{\'e}sentations automorphes
cohomologiques de $\GL_n(\Ade)$. De plus, si $\pi$ est une repr{\'e}sentation
automorphe cuspidale cohomologique de $\GL_n(\Ade)$ sur laquelle $\Ar$
agit trivialement, et si on peut prendre $W=\ungras$ dans la d{\'e}finition
ci-dessus, alors la repr{\'e}sentation de $\Hf$ associ{\'e}e {\`a} $\pi$ appara{\^i}t comme
un sous-quotient de $\H^*_c$.\footnote{On obtiendrait les repr{\'e}sentations
cohomologiques pour les $W$ non triviales en prenant la cohomologie {\`a}
coefficients dans un syst{\`e}me local non trivial sur les $X_K$.}

\end{theo}

\subsection{Cohomologie des vari{\'e}t{\'e}s de Shimura et conjecture \ref{conj:LC}}\label{1.4}

Les espaces localement sym{\'e}triques associ{\'e}s au groupe $\GL_n$ sont seulement
des vari{\'e}t{\'e}s analyiques r{\'e}elles, mais, en utilisant d'autres groupes, on
peut obtenir des espaces avec plus de structure (c'est-{\`a}-dire des vari{\'e}t{\'e}s
alg{\'e}briques sur des corps de nombres, appel{\'e}es \emph{vari{\'e}t{\'e}s de Shimura}).
Soit $\G$ un groupe
alg{\'e}brique r{\'e}ductif connexe sur $\Q$. \`A l'exception du th\'eor\`eme
de multiplicit\'e un fort,
toutes les d{\'e}finitions et les r{\'e}sultats
ci-dessus restent valables pour $\G$ 
(il faut remplacer $K_\infty=\SO(n)$ par un sous-groupe
compact connexe maximal de $\G(\R)$ et $\Ar$ par $\Se(\R)^\circ$, o{\`u}
$\Se$ est le sous-tore d{\'e}ploy{\'e} (sur $\Q$) maximal du centre de $\G$; pour
l'isomorphisme de Satake, il faut supposer $\G$ non ramifi{\'e} en $p$ et
remplacer $\GL_n(\Z_p)$ par un sous-groupe compact maximal hypersp{\'e}cial
de $\G(\Q_p)$, voir la section 1.10 de l'article \cite{T} de Tits).

Si par exemple $\G$ est le groupe symplectique $\Sp_{2n}\subset\GL_{2n}$
de la forme symplectique $x_1 y_{2n}+\dots x_n y_{n+1}-x_{n+1}y_n-\dots
-x_1 y_{2n}$ et $K=\Ker(\G(\widehat{\Z})\fl\G(\Z/N\Z))$ pour $N$ un entier
$\geq 3$ (pour que $K$ soit assez petit), alors l'espace localement sym{\'e}trique
associ{\'e} $X_K^{\G}$ est l'espace de modules des vari{\'e}t{\'e}s ab{\'e}liennes de dimension
$n$ sur $\C$ principalement polaris{\'e}es et munies d'une structure de niveau
$N$ (voir la section 11 de l'article \cite{K1} de Kottwitz). On peut d{\'e}finir
le probl{\`e}me de modules sur $\Q$ (ou m{\^e}me $\Z$), et Mumford a montr{\'e} que
ce probl{\`e}me de modules est repr{\'e}sentable par un sch{\'e}ma quasi-projectif
(cf. le th{\'e}or{\`e}me 7.9 de \cite{GIT}). Les correspondances de Hecke ont
aussi une description modulaire, et sont donc d{\'e}finies sur $\Z$.
En utilisant le th{\'e}or{\`e}me de comparaison
entre cohomologie de Betti et cohomologie {\'e}tale et en utilisant
$\oQ_\ell$ au lieu de $\C$ comme corps de coefficients, on en d{\'e}duit que
$H^*_{c,\G}:=\varinjlim_{K}H^*_c(X_K^{\G})$
est muni d'une action de $\Gal(\oQ/\Q)$ qui
commute {\`a} l'action de l'alg{\`e}bre de Hecke $\Hf$. On peut donc {\'e}crire la
semi-simplifi{\'e}e de cette repr{\'e}sentation de la mani{\`e}re suivante :
\[(H^i_{c,\G})^{ss}=\bigoplus_{\pi=\pi_\infty\otimes\pi_f}\pi_f\otimes\sigma^i(\pi_f),\]
o{\`u} $\pi$ parcourt l'ensemble des repr{\'e}sentations
automorphes cohomologiques de $\G(\Ade)$ et les $\sigma^i(\pi_f)$ sont des
repr{\'e}sentations semi-simples de $\Gal(\oQ/\Q)$. Il n'est pas {\'e}vident en
g{\'e}n{\'e}ral de d{\'e}terminer $\sigma^i(\pi_f)$ en fonction de $\pi_f$, mais, si
l'on utilise la cohomologie d'intersection (voir par exemple
\cite{BBD}) au lieu de la cohomologie {\`a} supports compacts, alors on une
formule conjecturale tr{\`e}s pr{\'e}cise pour $\sigma^i(\pi_f)$, due {\`a} Langlands,
Rapoport et Kottwitz (cf. la section 10 de \cite{K1}). Cette formule fait
intervenir la conjecture de r{\'e}ciprocit{\'e} de Langlands pour le groupe $\G$.

De plus, tout le paragraphe pr{\'e}c{\'e}dent
est en fait valable pour les groupes $\G$ dont les espaces
localement sym{\'e}triques sont des vari{\'e}t{\'e}s de Shimura de type PEL, par
exemple les groupes unitaires et certains groupes orthogonaux (voir la
section 5 de l'article \cite{K2} de Kottwitz), {\`a} la diff{\'e}rence que la structure
de vari{\'e}t{\'e} alg{\'e}brique des $X_K^{\G}$ est en g{\'e}n{\'e}ral d{\'e}finie non pas sur $\Q$,
mais sur une extension finie de $\Q$ appel{\'e}e corps reflex. 
Si on choisit le groupe $\G$ (et le degr{\'e} $i$) correctement, la cohomologie
d'intersection en degr{\'e} $i$ (qui, dans le cas o{\`u} $\G^{der}$ est anisotrope
sur $\Q$, est simplement $\H^i_{c,\G}$, car les $X_K^{\G}$ sont alors des
sch{\'e}mas projectifs) r{\'e}alise conjecturalement une partie de
la correspondance de Langlands
pour le groupe $\G$. De plus, dans de nombreux cas, la conjecture est en fait
connue.

Si l'on veut obtenir des informations sur la conjecture \ref{conj:LC} pour
$\GL_n$, on peut passer des repr{\'e}sentations automorphes de $\GL_n$ {\`a} celles
d'un autre groupe en utilisant le principe de fonctorialit{\'e} de Langlands
(conjectural lui aussi en g{\'e}n{\'e}ral, mais connu dans les cas que l'on veut
utiliser ici). Soit $\G$ un groupe r{\'e}ductif connexe sur $\Q$, d{\'e}ploy{\'e} pour
simplifier, et soit
$\widehat{\G}$ son dual de Langlands (le groupe r{\'e}ductif connexe sur $\C$
obtenu en {\'e}changeant le r{\^o}le des racines et des coracines dans
la donn{\'e}e radicielle de $\G$). Alors les repr{\'e}sentations automorphes cuspidales
de $\G(\Ade)$ sont conjecturalement classifi{\'e}es par les param{\`e}tres de
Langlands, qui sont des morphismes alg{\'e}briques irr{\'e}ductibles\footnote{
C'est-{\`a}-dire dont l'image n'est contenue dans aucun sous-groupe parabolique.}
$\Lf_\Q\fl\widehat{\G}$. En g{\'e}n{\'e}ral, cette param{\'e}trisation n'est plus bijective,
et l'on s'attend {\`a} ce que chaque param{\`e}tre corresponde {\`a} un ensemble fini
de repr{\'e}sentations automorphes cuspidales, appel{\'e} un $L$-paquet.
En tout cas, si l'on a deux groupes r{\'e}ductifs connexes~$\G$ et $\GH$ et un
morphisme $\widehat{\G}\fl\widehat{\GH}$, cette conjecture implique que l'on
a un {\og transfert\fg} qui envoie une repr{\'e}sentation automorphe cuspidale 
de $\G(\Ade)$ sur un $L$-paquet de repr{\'e}sentations automorphes
de $\GH(\Ade)$ (de
mani{\`e}re compatible {\`a} l'isomorphisme de Satake aux places o{\`u} les
repr{\'e}sentations sont non ramifi{\'e}es). On devrait aussi pouvoir caract{\'e}riser
l'image de ce transfert.\footnote{Le transfert
ne pr{\'e}serve pas la cuspidalit{\'e} en g{\'e}n{\'e}ral, il faut donc travailler avec
les repr{\'e}sentations automorphes discr{\`e}tes, qui sont (conjecturalement)
param{\'e}tr{\'e}es par des
param{\`e}tres d'Arthur $\psi:\Lf_\Q\times\SL_2(\C)\fl\widehat{\G}$ et non des
param{\`e}tres de Langlands.}

Par exemple, si $\G=\Sp_{2n}$ et $\GH=\GL_{2n+1}$, alors $\widehat{\G}=\SO_{2n+1}
(\C)$ se plonge de mani{\`e}re {\'e}vidente dans $\widehat{\GH}=\GL_{2n+1}(\C)$. Dans
ce cas, l'existence du transfert et ses propri{\'e}t{\'e}s ont {\'e}t{\'e} {\'e}tablies par
Arthur dans le livre \cite{A}; l'image du transfert est caract{\'e}ris{\'e}e par une
condition d'autodualit{\'e} et une condition sur les p{\^o}les d'une certaine
fonction~$L$, voir le th{\'e}or{\`e}me 1.5.3 de \cite{A}. Si $\G$ est un
groupe unitaire et $\GH$ un groupe g\'en\'eral lin\'eaire, on a des r\'esultats
similaires, dus \`a Mok (\cite{M}) et Kaletha-Minguez-Shin-White (\cite{KMSW}).

En utilisant le transfert des groupes unitaires vers les groupes g{\'e}n{\'e}raux
lin{\'e}aires, le calcul de la cohomologie des vari{\'e}t{\'e}s de Shimura de certains
groupes unitaires, le lemme fondamental et des techniques d'interpolation
$p$-adique pour attraper certaines repr{\'e}sentations automorphes, on arrive
au r{\'e}sultat suivant (d\^u, au moins, {\`a} Kottwitz, Clozel, Labesse, Harris-Taylor,
Fargues, Mantovan, Shin, Laumon-Ng\^o, Waldspurger,
Bella{\"i}che-Chenevier\footnote{L'auteur
regrette de ne pouvoir garantir l'exhaustivit{\'e} de cette liste.}) :

\begin{theo}[\cite{PL}]\footnote{On a un r{\'e}sultat similaire pour les
repr{\'e}sentations du groupe $\GL_n(\Ade_F)$, o{\`u} $F$ est un corps de nombres
totalement r{\'e}el ou CM. Un {\'e}nonc{\'e} pr{\'e}cis est rappel{\'e} dans le th{\'e}ror{\`e}me V.1.4
de \cite{S}.}
La conjecture \ref{conj:LC} est vraie pour les repr{\'e}sentations
automorphes cuspidales cohomologiques autoduales (c'est-{\`a}-dire isomorphes
{\`a} leur contragr{\'e}diente).
\end{theo}

\subsection{Repr{\'e}sentations non autoduales}
\label{1.5}

Les m{\'e}thodes de la section pr{\'e}c{\'e}dentes ne peuvent s'appliquer aux repr{\'e}sentations
non autoduales de $\GL_n(\Ade)$. En effet, toutes les repr{\'e}sentations venant
par transfert depuis un groupe ayant une vari{\'e}t{\'e} de Shimura v{\'e}rifient une
propri{\'e}t{\'e} d'autodualit{\'e}.

L'approche suivante a {\'e}t{\'e} sugg{\'e}r{\'e}e par Clozel.
Si $\G=\Sp_{2n}\subset\GL_{2n}$ est
le groupe symplectique de la forme antidiagonale (comme plus haut), alors
le groupe \mbox{$\Pa=\Sp_{2n}\cap\begin{pmatrix}* & * \\ 0 & *\end{pmatrix}$} est un
sous-groupe parabolique maximal de $\G$, de quotient de Levi isomorphe
{\`a} $\GL_n$. Si $\pi$ est une repr{\'e}sentation automorphe cuspidale cohomologique
de $\GL_n(\Ade)$, elle d{\'e}finit donc par induction parabolique une repr{\'e}sentation
automorphe\footnote{Ou plusieurs...} $\Pi$
de $\Sp_{2n}(\Ade)$, cf. \cite{L}, qui n'est {\'e}videmment pas cuspidale, mais se
trouve {\^e}tre aussi cohomologique. Comme le groupe $\Sp_{2n}$ admet une
vari{\'e}t{\'e} de Shimura, on peut essayer d'appliquer les techniques de la section
pr{\'e}c{\'e}dente {\`a} $\Pi$. Malheureusement, les repr{\'e}sentations galoisiennes
qui apparaissent, dont la composante $\Pi_f$-isotypique (o{\`u} $\Pi=
\Pi_\infty\otimes\Pi_f$) de la cohomologie
des $X_K^{\G}$, ne sont pas tr\`es int\'eressantes
(on obtient quelque chose qui ressemble beaucoup {\`a} la puissance ext{\'e}rieure
$n$-i{\`e}me de la repr{\'e}sentation de dimension~$n$ que l'on essaie de construire).

Une autre id{\'e}e, au lieu d'utiliser directement la cohomologie de Betti
des $X_K^{\G}$, est d'utiliser l'autre r{\'e}alisation cohomologique des formes
automorphes (comme sections de certains fibr{\'e}s vectoriels, dits {\og automorphes\fg},
sur les espaces localement sym{\'e}triques) pour approcher $p$-adiquement la
repr{\'e}sentation
$\Pi$ par des repr{\'e}sentations automorphes cuspidales cohomologiques de
$\Sp_{2n}(\Ade)$ (o\`u $p$ est un nombre premier fix\'e arbitrairement).
Ces repr{\'e}sentations se transf{\`e}rent alors en des
repr{\'e}sentations automorphes autoduales de $\GL_{2n+1}(\Ade)$ (voir la
discussion sur le transfert plus haut), qui ont des repr{\'e}sentations
galoisiennes associ{\'e}es ({\`a} coefficients dans $\oQ_p$). En prenant la limite
de ces repr{\'e}sentations galoisiennes (ou plut{\^o}t de leurs caract{\`e}res), on
obtient une repr{\'e}sentation galoisienne de dimension $2n+1$, dont il est
possible d'extraire la repr{\'e}sentation $\rho_\pi$ cherch{\'e}e. Cette strat{\'e}gie
a {\'e}t{\'e} men{\'e}e {\`a} bien de mani{\`e}re ind{\'e}pendante par Harris-Lan-Taylor-Thorne 
(\cite{HLTT}),
Scholze (\cite{S}) et Boxer. En fait, Scholze et Boxer prouvent des
r{\'e}sultats plus forts, que nous allons expliquer ci-dessous.

D'abord, remarquons que le paragraphe pr{\'e}c{\'e}dent n'a pas de sens a priori,
car les repr{\'e}sentations automorphes sont {\`a} coefficients dans $\C$, et
non $\oQ_p$. Mais d'apr{\`e}s le th{\'e}or{\`e}me \ref{th:Franke} ci-dessus (qui est valable
pour un groupe r{\'e}ductif connexe quelconque), on peut
remplacer l'{\'e}tude des repr{\'e}sentations automorphes cohomologiques par celles
des sous-quotients de la repr{\'e}sentation $H^*_{c,\G}=\varinjlim_{K}\H^*_c(X_K^{\G})$
de l'alg{\`e}bre de Hecke globale $\Hf_{\G}$ de $\G$ (l'alg{\`e}bre des fonctions
localement constantes {\`a} support compact $\G(\Af)\fl\C$, munie du produit
de convolution). Mais tant la cohomologe de Betti que l'alg{\`e}bre de Hecke
globale ont un sens si l'on remplace le corps de coefficients $\C$ par un
anneau commutatif de coefficients quelconque $A$;\footnote{Ce n'est pas
tout \`a fait vrai. A priori, il
faudrait avoir une mesure de Haar sur $\G(\Af)$ \`a valeurs dans $A$, ce qui
est une condition non triviale. On verra dans le paragraphe suivant comment
faire si $A=\Z$.} 
notons $\H^*_{c,\G}(A)$
et $\Hf_{\G}(A)$ les objets ainsi obtenus. Si l'on a deux repr\'esentations
de $\Hf_{\G}(\Q_p)$ apparaissant dans $\H^*_{c,\G}(\Q_p)$, cela a un sens
de demander qu'elles soient $p$-adiquement proches, mais ce n'est pas encore
exactement ce que l'on veut faire. En effet, on veut approximer la classe
d'isomorphisme d'une repr\'esentation et non la repr\'esentation elle-m\^eme.
Il est donc naturel de chercher \`a approximer le caract\`ere de la
repr\'esentation, et pour cela il est plus commode de fixer le niveau $K$
(afin d'avoir des repr\'esentations de dimension finie).

On fixe donc un sous-groupe compact ouvert $K$ de $\G(\Af)$, et on suppose
que $K$ est de la forme $\prod_v K_v$, avec $v$ parcourant les nombres premiers
et $K_v$ un sous-groupe compact ouvert de $\G(\Q_v)$.
Il existe un
ensemble fini $S$ de nombres premiers tel que, pour tout $v\not\in S$,
le groupe $\G$ soit non ramifi\'e en $v$ et $K_v$ soit hypersp\'ecial (voir
le d\'ebut de la section \ref{1.4} pour la d\'efinition), et on fixe un tel $S$. Pour
tout $v\not\in S$, on note $\Hf_{v,\G}$ l'alg\`ebre de
Hecke locale non ramifi\'ee en $v$, c'est-\`a-dire l'alg\`ebre des
fonctions $f:\G(\Q_v)\fl\Z$ \`a support compact et bi-invariantes par
$K_v$, munie du produit de convolution (pour la mesure de Haar sur
$\G(\Q_v)$ telle que $K_v$ soit de volume $1$; pour le fait que le
produit de convolution de deux fonctions \`a valeurs dans $\Z$ est bien une
fonction \`a valeurs dans $\Z$, voir par exemple la section 2 de \cite{G}).
On note aussi $\Hf^S_{\G}=\bigotimes_{v\not\in S}\Hf_{v,\G}$. Gr\^ace \`a
l'isomorphisme de Satake, les alg\`ebres $\Hf_{v,\G}$ et $\Hf^S_{\G}$ sont
commutatives, donc leurs repr\'esentations irr\'eductibles sont simplement
des caract\`eres.
De plus, l'alg\`ebre $\Hf^S_{\G}$ agit sur la cohomologie (de Betti) $\H^*_c(X_K^{\G},A)$, pour tout anneau commutatif $A$.
Notons que ces constructions sont possibles pour n'importe quel groupe
r\'eductif connexe $\G$.

On peut alors faire la chose suivante. On prend comme avant $\G=\Sp_{2n}$,
et on voit $\GL_n$ comme le quotient de Levi d'un sous-groupe parabolique
maximal de $\G$.
Soient~$\pi$ et $\Pi$ comme plus haut.
On choisit le sous-groupe compact ouvert $K$ de $\G(\Af)$ tel que $\Pi^\K\not=0$.
Alors $\Pi$ correspond \`a un caract\`ere $\varphi$ de $\Hf^S_{\G}$ qui appara\^it comme
un sous-quotient de $\H^*_c(X_K^{\G},\C)$. Le choix d'un isomorphisme $\C\simeq
\oQ_p$ permet de voir $\varphi$ comme un caract\`ere \`a valeurs dans
$\oQ_p$, et on montre qu'il existe une extension finie $E$ de $\Q_p$ telle que
$\varphi$ soit \`a valeurs dans $\Of_E$. La question devient alors de savoir
si l'on peut
trouver une suite $(\Pi_i)_{i\in\Nat}$
de repr\'esentations automorphes cohomologiques \emph{cuspidales} de
$\G(\Ade)$ telle que $(\Pi_i)^{K^S}\not=0$ pour tout $i$ (o\`u
$K^S=\prod_{v\not\in S}K_v$) et que, si
$\varphi_i:\Hf^S_{\G}\fl\oQ_p$ est le caract\`ere associ\'e \`a $\Pi_i$ comme
plus haut, on ait
\[\lim_{i\rightarrow\infty}\varphi_i(x)=\varphi(x)\]
pour tout $x\in\Hf^S_{\G}$, o\`u on prend la limite pour la topologie $p$-adique sur
$\oQ_p$.

Harris, Lan, Taylor et Thorne ont \'et\'e les premiers \`a donner une
r\'eponse affirmative \`a cette question, ce qui leur a permis de prouver
le th\'eor\`eme suivant :

\begin{theo}[\cite{HLTT}]\label{th:HLTT}
La conjecture \ref{conj:LC} est vraie si
$\pi$ est cohomologique.\footnote{Et on a un r\'esultat similaire pour les
repr\'esentations automorphes de $\GL_n(\Ade_F)$, si $F$ est un corps de
nombres totalement r\'eel ou CM. Notons que la preuve dans le cas CM utilise les
vari\'et\'es de Shimura des groupes unitaires quasi-d\'eploy\'es au lieu
de celles des groupes symplectiques.}
\end{theo}

Rappelons que la m\'ethode esquiss\'ee ci-dessus ne donne pas directement la
repr\'esentation $\rho_\pi$, mais plut\^ot quelque chose qui ressemble
\`a $\rho_\pi\oplus(\rho_\pi)^*$. Il y a une derni\`ere \'etape qui consiste
\`a extraire $\rho_\pi$ de cette repr\'esentation, et que nous ignorerons
totalement (voir la section V.3 \cite{S}).

Scholze et Boxer ont reprouv\'e ce r\'esultat (ind\'ependamment de Harris-Lan-Taylor-Thorne
et ind\'ependamment l'un de l'autre) et l'ont g\'en\'eralis\'e aux caract\`eres
de $\Hf^S_{\GL_n}$ apparaissant dans la torsion de $H^*_c(X_K^{\GL_n},\Z)$. (Noter
qu'il s'agit maintenant de l'alg\`ebre de Hecke et de l'espace localement
sym\'etrique pour $\GL_n$, et non $\Sp_{2n}$.) Plus pr\'ecis\'ement, on a la
conjecture suivante, due \`a Ash :

\begin{conj}[\cite{As1},\cite{A2}] \label{conj:A}Soit $S$ un ensemble fini
de nombres premiers. Si $\varphi:\Hf^S_{\GL_n}\fl\Fi_p$ est un caract\`ere
qui appara\^it dans $\H^*_c(X_K^{\GL_n},\Fi_p)$, pour $K=\prod_v K_v$ un sous-groupe
compact ouvert de $\GL_n(\Af)$ tel que $K_v=\GL_n(\Z_v)$ si $v\not\in S$,
alors il existe une repr\'esentation semi-simple
$\rho:\Gal(\oQ/\Q)\fl\GL_n(\Fi_p)$ telle que, pour tout $v\not\in S\cup\{p\}$,
$\varphi_{|\Hf_{v,\GL_n}}$ et $\rho_{|\Gal(\oQ_p/\Q_p)}$ se correspondent par
l'isomorphisme de Satake.\footnote{Voir les explications apr\`es la conjecture
\ref{conj:LC}.}

\end{conj}

\section{\'Enonc\'e du th\'eor\`eme principal et strat\'egie de la preuve}

\begin{theo}[Scholze \cite{S}, Boxer]\label{th}
 La conjecture \ref{conj:A} est
vraie.\footnote{Comme dans le cas du th\'eor\`eme \ref{th:HLTT}, Scholze et
Boxer prouvent en fait un r\'esultat valable pour le groupe $\GL_n$ sur
un corps de nombres totalement r\'eel ou CM.}

\end{theo}

Notons que l'on a en fait une version du th\'eor\`eme ci-dessus pour la
cohomologie \`a coefficients dans $\Z/p^m\Z$, pour tout entier strictement
positif $m$ (mais l'\'enonc\'e est un peu plus compliqu\'e \`a formuler, voir le
th\'eor\`eme V.3.1 de \cite{S}). En passant \`a la limite sur $m$, on obtient
donc une nouvelle preuve du th\'eor\`eme \ref{th:HLTT}. Dans ce texte, nous
nous concentrerons pour simplifier sur le cas o\`u $m=1$.

On va pr\'esenter la preuve de Scholze (dans les grandes lignes),
qui utilise la th\'eorie des
espaces perfecto\"ides (voir l'article introductif \cite{S2} de Scholze ou
l'expos\'e \cite{Fo} de Fontaine au s\'eminaire Bourbaki). La preuve de
Boxer n'utilise pas cette th\'eorie, mais les d\'etails de cette preuve ne
sont pas connus de l'auteur.

On se ram\`ene tout d'abord \`a un \'enonc\'e sur la torsion dans la
cohomologie de certaines vari\'et\'es de Shimura, de la mani\`ere suivante.
Comme dans la section \ref{1.5}, on pose \mbox{$\G=\Sp_{2n}$,} et on voit
$\GL_n$ comme le quotient de Levi d'un sous-groupe parabolique maximal $\Pa$
de $\G$. Pour tout nombre premier $v$, l'application {\og terme constant le
long de $\Pa$\fg}
(voir par exemple la formule (19) de \cite{C}) d\'efinit un morphisme
injectif d'alg\`ebres $\Hf_{v,\G}\fl\Hf_{v,\GL_n}$. D'o\`u, si $S$ est un
ensemble fini de nombres premiers, un morphisme $\Hf^S_{\G}\fl
\Hf^S_{\GL_n}$.

D'autre part, soit $K$ un sous-groupe ouvert compact de $\G(\Af)$. L'espace
localement sym\'etrique $X_K^{\G}$ n'est pas compact, mais il admet
une compactification $\overline{X}_K^{\G}$ appel\'ee \emph{compactification
de Borel-Serre} et d\'efinie dans \cite{BS}, qui est une vari\'et\'e
analytique r\'eelle \`a coins ayant la m\^eme cohomologie
(sans supports) que $X_K^{\G}$
et telle que le bord $\overline{X}_K^{\G}-X_K^{\G}$ admette
une stratification par des sous-vari\'et\'es analytiques r\'eelles de la
forme $X_{K_Q}^{\QP}$, pour $\QP$ un sous-groupe parabolique de $\G$ et $K_Q$
un sous-groupe compact ouvert de $\QP(\Af)$. En particulier, on a
des strates correspondant au sous-groupe parabolique maximal $\Pa$.
Comme on a un morphisme surjectif $\pi:\Pa\fl\GL_n$ (qui identifie $\GL_n$
au quotient de Levi de $\Pa$), on obtient, pour tout sous-groupe compact
ouvert $K_P$ de $\Pa(\Af)$, un morphisme surjectif $X_{K_P}^{\Pa}\fl
X_{\pi(K_P)}^{\GL_n}$, qui se trouve \^etre un fibr\'e en $(S^1)^N$, o\`u $N$
est la dimension du radical unipotent de $\Pa$. En utilisant ce fait
et la suite exacte longue d'excision, on montre que tout caract\`ere
$\varphi$ de $\Hf^S_{\GL_n}$ qui appara\^it dans un $H^*_c(X_{K_{\GL_n}}^{\GL_n},
\Fi_p)$ appara\^it aussi dans un $H^*_c(X_K^{\G},\Fi_p)$ (c'est-\`a-dire
que le compos\'e $\varphi':\Hf^S_{\G}\fl\Fi_p$ de $\varphi$ et du
morphisme $\Hf^S_{\G}\fl\Hf^S_{\GL_n}$ ci-dessus
appara\^it dans
un $H^*_c(X_K^{\G},\Fi_p)$). Pour les d\'etails, voir le d\'ebut de la section
V.2 de \cite{S}.

On se ram\`ene donc \`a montrer le th\'eor\`eme suivant :

\begin{theo}[Th\'eor\`emes I.5 et IV.3.1 de \cite{S}]\label{th:princ}\footnote{Ce th\'eor\`eme 
est en fait valable pour tous les groupes $\G$ d\'efinissant des vari\'et\'es 
de Shimura de type Hodge.}
Soit $\varphi:\Hf^K_{\G}\fl\Fi_p$ un caract\`ere apparaissant dans un
$H^*_c(X_K^{\G},\Fi_p)$ (pour un $K$ hypersp\'ecial aux places hors de~$S$,
comme au-dessus du th\'eor\`eme \ref{th:HLTT}). Alors il existe
une repr\'esentation automorphe cuspidale cohomologique $\pi$ de
$\G(\Ade)$, non ramifi\'ee en les places hors de $S$ et telle que, si
$\psi$ est le
caract\`ere correspondant de $\Hf^S_{\G}$, vu comme un morphisme
$\Hf^S_{\G}\fl\overline{\Z}_p$, alors la r\'eduction modulo $p$ de $\psi$
est \'egale \`a $\varphi$.

\end{theo}

Pour prouver ce th\'eor\`eme, on veut
comparer un caract\`ere de $\Hf^S_{\G}$ qui appara\^it dans
un groupe de cohomologie de Betti (ou \'etale) $H^*_c(X_K^{\G},\Fi_p)$
avec un caract\`ere venant d'une repr\'esentation automorphe cuspidale.
On peut voir les formes automorphes cuspidales sur $\G(\Ade)$ comme les
sections d'un certain faisceau coh\'erent sur $X_K^{\G}$, et Scholze
a justement un th\'eor\`eme de comparaison entre la cohomologie d'un
$\Fi_p$-syst\`eme local $\mathbb{L}$ sur un espace adique
propre et lisse $X$ sur $\C_p$
et celle de $\mathbb{L}\otimes\Of^+_X/p$ (corollaire 5.11 de
\cite{S1} et th\'eor\`eme 3.3 de \cite{S2}). Voir la section \ref{comp}
pour des rappels sur ce th\'eor\`eme.
Dans notre cas, le syst\`eme local $\mathbb{L}$ sera le syst\`eme local
trivial, et on compare sa cohomologie \`a celle du  {\og faisceau des formes
cuspidales\fg} sur $X_K^{\G}$, ou plut\^ot sur sa compactification de
Baily-Borel $X_{K}^{\G,*}$. Voir la section \ref{cc}.

Il faut ensuite passer de la cohomologie du faisceau des formes cuspidales
\`a ses sections globales. Or, Scholze a prouv\'e que la cohomologie d'un
faisceau coh\'erent sur un espace perfecto\"ide affino\"ide est presque
nulle. Les vari\'et\'es de Shimura ne sont pas perfecto\"ides, mais
Scholze prouve le r\'esultat suivant : Soit $p$ un nombre
premier. Si on fixe
un sous-groupe compact ouvert $K^p$ de $\G(\Af^p)$, o\`u $\Af^p=\prod'_{v\not=p}
\Q_v$, alors la limite projective 
$X_{K^p}^{\G}$ des $X_{K_p K^p}^{\G}$ lorsque $K_p$ parcourt
les sous-groupes compacts ouverts de $\G(\Q_p)$  {\og est\fg} un espace
perfecto\"ide (dans un sens \`a pr\'eciser), qu'on appelle parfois  {\og vari\'et\'e
de Shimura de niveau infini en $p$\fg}.
De plus,
on a un r\'esultat similaire
pour les compactifications de Baily-Borel de ces vari\'et\'es. 
La preuve de ce r\'esultat utilise de mani\`ere essentielle le morphisme
des p\'eriodes de Hodge-Tate, qui n'est d\'efini que sur la vari\'et\'e
de Shimura de niveau infini. Voir la section \ref{HT}.

On se place donc sur la vari\'et\'e de Shimura de niveau infini en $p$.
En utilisant un recouvrement (explicite) par des ouverts affino\"ides
perfecto\"ides, on montre que l'on peut approximer la cohomologie du
faisceau des formes cuspidales par des sections de ce faisceau sur ces
ouverts, qui sont des vecteurs propres pour l'action de l'alg\`ebre de
Hecke $\Hf^S_{\G}$.
Il faut encore prolonger ces sections \`a toute la vari\'et\'e
de Shimura, sans changer les valeurs propres pour l'action de l'alg\`ebre de
Hecke. La m\'ethode classique utilise l'invariant de Hasse, qui ne suffit
pas ici. Cependant, en utilisant le morphisme des p\'eriodes de
Hodge-Tate (qui est \'equivariant sous l'action des correspondances de Hecke
en dehors de $p$), on construit de  {\og faux\fg} invariants de Hasse qui
jouent le m\^eme r\^ole et permettent de finir la preuve. Voir la section
\ref{faux}.

Notons que la preuve du th\'eor\`eme \ref{th} n'est pas l'unique
application des m\'ethodes de \cite{S}. Voir la section \ref{appl} pour
quelques autres applications.

\section{Un th\'eor\`eme de comparaison}\label{comp}

Le th\'eor\`eme suivant est le d\'ebut de la preuve par Scholze du th\'eor\`eme
de comparaison entre cohomologie de de Rham et cohomologie \'etale $p$-adique.

\begin{theo}[Th\'eor\`eme 5.1 de \cite{S1} ou th\'eor\`eme 3.3 de \cite{S2}]
\label{th:comp1}
Soit $C$ une extension compl\`ete alg\'ebriquement close de $\Q_p$,
d'anneau des entiers $\Of_C$, et soit $X$ un espace adique propre et lisse
\footnote{En fait, l'hypoth\`ese de lissit\'e est superflue, voir le th\'eor\`eme 3.17 de \cite{S2}.}
sur $(C,\Of_C)$. Alors, pour tout $\Fi_p$-syst\`eme local $\mathbb{L}$
sur $X$, on a des presque isomorphismes
\[H^i(X_{\et},\mathbb{L})\otimes\Of^a/p\simeq\H^i(X_{\et},\Of_X^{+a}/p),\]
o\`u $\Of_X^+\subset\Of_X$ est le faisceau des fonctions born\'ees par
$1$ sur $X$.

\end{theo}

Pour des rappels et des r\'ef\'erences sur les espaces adiques, voir le
d\'ebut de la section~2 de \cite{Fo}. On peut par exemple prendre pour $X$
l'espace adique associ\'e \`a un sch\'ema propre et lisse sur $C$, mais une
des forces du r\'esultat de Scholze est qu'il s'applique aussi aux espaces
adiques ne venant pas d'un sch\'ema. Pour des rappels sur le langage des
presque math\'ematiques, voir la section 1.6 de \cite{Fo}; un {\og presque
isomorphisme\fg} (ou un isomorphisme de presque $\Of_C$-modules)
est un morphisme de $\Of_C$-modules dont le noyau et le
conoyau sont annul\'es par l'id\'eal maximal de $\Of_C$, et les {\og $a$\fg} en
exposant sont l\`a pour rappeler que l'on consid\`ere les objets comme des
presque $\Of_C$-modules.

Mentionnons aussi que le th\'eor\`eme ci-dessus est toujours valable si l'on
remplace $\Of_C$ par un sous-anneau de valuation ouvert et born\'e $C^+$ de
$C$, et dans le cas relatif (c'est-\`a-dire pour les images directes par
un morphisme propre et lisse de vari\'et\'es analytiques rigides sur $C$,
cf. le corollaire 5.11 de \cite{S1}).

\noindent{\sc Preuve} (esquisse) --- Il y a trois \'etapes dans la preuve.
\begin{enumerate}
\item[(1)] (Cf. le th\'eor\`eme 4.9 de \cite{S1}.)
Si $X$ est une vari\'et\'e analytique rigide affino\"ide et
connexe sur $C$ (o\`u $C$ est comme dans l'\'enonc\'e), si $x\in X(C)$ et si
$\pi=\pi_1(X,x)$ (le groupe fondamental \'etale profini de $X$), alors, pour
tout $\Fi_p$-syst\`eme local $\mathbb{L}$ sur~$X$, les morphismes canoniques
\[H^i_{cont}(\pi,\mathbb{L}_x)\fl H^i(x_{\et},\mathbb{L})\]
sont des isomorphismes. (Autrement dit, $X$ est un $K(\pi,1)$ pour les coefficients
de $p$-torsion.)

Il faut montrer que toute classe de cohomologie dans un $H^i(X_{\et},\mathbb{L})$
pour $i>0$ est tu\'ee par un rev\^etement fini \'etale de $X$. On se ram\`ene
facilement au cas o\`u $\mathbb{L}$~est le syst\`eme local trivial $\Fi_p$.

Chaque rev\^etement fini \'etale de $X$ est affino\"ide, donc de la forme
$\Spa(A,A^+)$. En prenant une compl\'etion appropri\'ee de
la limite inductive de ces $C$-alg\`ebres~$A$, on obtient une
$C$-alg\`ebre perfecto\"ide $(A_\infty,A_\infty^+)$, dont le $\Spa$ est
un espace perfecto\"ide $X_\infty$, qui m\'erite le nom de rev\^etement
fondamental de $X$. (Voir \cite{Fo} pour une introduction aux espaces
perfecto\"ides.) Il s'agit maintenant de montrer que $\H^i(X_{\infty,\et},
\Fi_p)=0$ pour $i>0$. En utilisant le basculement (ou tilt), on se ram\`ene
\`a montrer le r\'esultat similaire pour l'espace perfecto\"ide $X_\infty^\flat$
sur le corps perfecto\"ide~$C^\flat$, qui est de caract\'eristique $p$.
Cela r\'esulte alors de la suite exacte d'Artin-Schreier $0\rightarrow\Fi_p
\rightarrow\Of_{X_\infty^\flat}\rightarrow\Of_{X_\infty^\flat}\rightarrow 0$ et
du fait que $X_\infty^\flat$ n'a pas de rev\^etement fini \'etale non trivial.

\item[(2)] On prouve ensuite que, si $X$ est un espace adique propre et
lisse sur $C$ et $\mathbb{L}$ est un $\Fi_p$-syst\`eme local sur $X$, alors
les groupes de cohomologie $H^i(X_{\et},\mathbb{L}\otimes\Of_X^+/p)$
sont presque de type fini, et presque nuls pour $i>2\dim X$ (lemme 5.8 de
\cite{S1}). L'id\'ee de la preuve est classique (et d\'ej\`a utilis\'ee par
Cartan et Serre pour les vari\'et\'es analyiques complexes et par Kiehl pour
les vari\'et\'es analytiques rigides; bien s\^ur, il faut faire marcher cette
id\'ee) : on calcule la cohomologie de $X$ en utilisant le complexe de {\v C}ech
d'un nombre fini de recouvrements par des ouverts affino\"ides dont chacun
raffine assez le pr\'ec\'edent. (Les m\'ethodes de (1) sont utilis\'ees pour
prouver que tous les groupes de cohomologie qui apparaissent sont
presque de type fini.)

\item[(3)] Enfin, pour prouver le th\'eor\`eme, on utilise la topologie
pro-\'etale de $X$. C'est une topologie de Grothendieck plus fine que la
topologie \'etale; l'id\'ee de base est simplement
que l'on autorise des recouvrements par des
limites projectives d'espaces \'etales sur $X$, mais les d\'etails techniques
ne sont pas totalement \'evidents (voir la section 3 de \cite{S1} pour la
d\'efinition rigoureuse). Ce qui rend cette topologie si utile est le fait que
\emph{tout espace adique localement noeth\'erien est pro-\'etale localement
perfecto\"ide} (voir la proposition 4.8 de \cite{S1}). On peut donc introduire
le faisceau structurel compl\'et\'e bascul\'e $\widehat{\Of}^+_{X^\flat}$, qui
est un faisceau (de $p$-torsion) sur le site pro-\'etale $X_{\proet}$, et
l'utiliser pour calculer les $H^i(X_{\proet},\mathbb{L})$ via la suite exacte
longue de cohomologie pro-\'etale de la suite exacte
d'Artin-Schreier $0\rightarrow\mathbb{L}\rightarrow\mathbb{L}\otimes
\widehat{\Of}^+_{X^\flat}\rightarrow\mathbb{L}\otimes\widehat{\Of}^+_{X^\flat}
\rightarrow 0$ (une suite exacte de faisceaux sur $X_{\proet}$). On utilise
le r\'esultat de finitude de (2) pour montrer que les morphismes
de connexion (c'est-\`a-dire ceux allant d'un $H^i$ dans un $H^{i+1}$) dans
la suite exacte longue ci-dessus sont nuls. Notons qu'on a aussi utilis\'e
sans le dire un th\'eor\`eme de comparaison entre cohomologies \'etale
et pro-\'etale, voir le corollaire 3.17 de \cite{S1}.
\end{enumerate}

Le th\'eor\`eme \ref{th:comp1} admet une version plus g\'en\'erale (au moins pour
les espaces adiques provenant de sch\'emas), pour des
syst\`emes de coefficients constructibles, qui est celle dont on aura besoin.

\begin{theo}[Th\'eor\`eme 3.13 de \cite{S2}]\label{th:comp2}
Soit $C$ comme
dans le th\'eor\`eme \ref{th:comp1}, soit $X$ l'espace adique associ\'e \`a
un sch\'ema propre sur $C$, et soit $\mathbb{L}$ l'image inverse
sur $X_{\et}$ d'un $\Fi_p$-faisceau \'etale constructible sur ce sch\'ema.
Alors on a des presque isomorphismes
\[H^i(X_{\et},\mathbb{L})\otimes\Of_C^a/p\simeq H^i(X_{\et},\mathbb{L}\otimes
\Of^{+a}_X/p).\]

\end{theo}

Comme pour le th\'eor\`eme \ref{th:comp1}, on peut remplacer $\Of_C$ par un
$C^+\subset C$ plus g\'en\'eral, et on a une version relative.

La preuve du th\'eor\`eme utilise, outre le th\'eor\`eme \ref{th:comp1} et la
r\'esolution des singularit\'es, le lemme simple suivant.

\begin{lemm}[\cite{S1}, lemme 3.14]
Soit $X$ un espace adique localement noeth\'erien sur
$\Spa(\Q_p,\Z_p)$, et soit $i:Z\fl X$ un sous-espace ferm\'e de $X$. Alors
le morphisme $i^*\Of_X^+/p\fl\Of_Z^+/p$ est un isomorphisme (de faisceaux sur
$X_{\et}$).

\end{lemm}

Pour l'application aux vari\'et\'es de Shimura, il est plus naturel d'introduire
d'abord la vari\'et\'e de Shimura perfecto\"ide.

\section{Vari{\'e}t{\'e} de Shimura perfecto{\"i}de et morphisme de Hodge-Tate}
\label{HT}

On utilise \`a nouveau les notations de la section \ref{1.4}, sauf que
l'on prend ici
\mbox{$\G=\GSp_{2n}\subset\GL_{2n}$} (le groupe g\'en\'eral symplectique).
Pour tout
sous-groupe ouvert compact $K$ de $\G(\Af)$ (assez petit), l'espace
localement sym\'etrique $X_K=X_K^{\G}$ est l'ensemble des points complexes d'une
vari\'et\'e quasi-projective lisse sur $\Q$, que l'on notera encore $X_K$.
Cette vari\'et\'e n'est pas projective (sauf si $n=0$), mais elle admet une
compactification canonique $X^*_K$, qui est une vari\'et\'e projective normale
sur $\Q$, appel\'ee compactification minimale, de Baily-Borel ou
de Satake-Baily-Borel (voir l'article \cite{BB} de Baily-Borel pour la
construction sur $\C$, et le livre \cite{CF} de Chai et Faltings pour la
construction sur $\Q$\footnote{Et m\^eme sur $\Z_p$ si $K=\G(\Z_p) K^p$ avec
$K^p\subset\G(\Af^p)$.}). Notons que $X^*_K$ est muni d'un
fibr\'e en droites ample canonique (sur $X_K$, c'est le d\'eterminant du
faisceau des formes diff\'erentielles 
invariantes de degr\'e~$1$ sur le sch\'ema
ab\'elien universel sur $X_K$).

On fixe un nombre premier $p$, et
on note $\X_K$ et $\X^*_K$ les espaces
adiques sur $\Spa(\Q_p,\Z_p)$ associ\'es aux sch\'emas $X_K$ et $X^*_K$.
On note encore $\omega$ le fibr\'e en droites sur $\X^*_K$ associ\'e au
fibr\'e en droites $\omega$ sur $X^*_K$.

D'autre part, soit $V=\Q^{2n}$, muni de la forme symplectique d\'efinie dans la
section \ref{1.4}. On note $\Flr$ la vari\'et\'e (sur $\Q$) des sous-espaces
totalement isotropes $W$ de dimension~$n$ de $V$, qui est munie d'un fibr\'e
en droites ample tautologique $\omega_\Flr=(\bigwedge^n W)^*$.
On note $\Fl$ l'espace adique
sur $\Spa(\Q_p,\Z_p)$ associ\'e \`a $\Flr$, et $\omega_\Fl$ le fibr\'e en droite
sur $\Fl$ correspondant \`a $\omega_\Flr$.

Enfin, on note $\Q_p^{\rm cycl}$ le compl\'et\'e de $\Q_p(\mu_{p^\infty})$ (qui est
l'extension de $\Q_p$ engendr\'ee par toutes les racines de $1$ d'ordre une
puissance de $p$), et $\Z_p^{\rm cycl}$ son anneau des entiers.
Notons que $\Q_p^{\rm cycl}$ est un corps perfecto\"ide.

Le th\'eor\`eme suivant est l'un des r\'esultats centraux de l'article
\cite{S}.

\begin{theo}[Th\'eor\`emes III.1.2 et III.3.17 de
\cite{S}]\label{th:perf}
 On fixe un sous-groupe compact ouvert assez petit $K^p$ de
$\G(\Af^p)$.

\begin{enumerate}
\item[\rm (i)] Il existe un espace perfecto\"ide $\X^*_{\Gamma(p^\infty)K^p}=
\X^*_{\Gamma(p^\infty)}$\footnote{Voir plus bas pour l'explication du {\og $\Gamma(
p^\infty)$\fg}.}
 sur $\Q_p^{cycl}$,
unique \`a isomorphisme unique pr\`es, tel que
\[\X^*_{\Gamma(p^\infty)K^p}\sim\varprojlim_{K_p}\X^*_{K_p K^p},\]
o\`u $K_p$ parcourt l'ensemble des sous-groupes compacts ouverts de $\G(\Q_p)$.
\footnote{Voir la d\'efinition 2.4.1 de \cite{SW} pour $\sim$. En particulier,
d'apr\`es le th\'eor\`eme 2.4.7 de \cite{SW},
ceci implique que le topos \'etale $\X^*_{\Gamma(p^\infty)K^p,\et}$ est la limite projective
des $\X^*_{K_p K^p,\et}$.}

\item[\rm (ii)] On a une application des p\'eriodes de Hodge-Tate
$\G(\Q_p)$-\'equivariante $\pi_{HT}:\X^*_{\Gamma(p^\infty)K^p}\fl\Fl$, qui commute avec
tous les op\'erateurs de Hecke hors de $p$,\footnote{C'est-\`a-dire donn\'e par
des $g\in\G(\Af)$ de composante en $p$ \'egale \`a $1$.} pour l'action triviale
de ces op\'erateurs sur $\Fl$.

\item[\rm (iii)] On a un isomorphisme canonique $\omega=\pi_{HT}^*\omega_\Fl$.

\item[\rm (iv)] On a un recouvrement de $\Fl$ par des ouverts affino\"ides $U$
\footnote{Explicites, ce sont les $\Fl_J$ d\'efinis plus bas.} tels que : 
\begin{enumerate}
\item $V=\pi_{HT}^{-1}(U)$ est perfecto\"ide affino\"ide;
\item pour tout $K_p\subset\G(\Q_p)$ assez petit, il existe $V_{K_p}\subset
\X^*_{K_p K^p}$ d'image inverse $V$ dans $\X^*_{\Gamma(p^\infty)K^p}$;
\item le morphisme suivant est d'image dense :
\[\varinjlim_{K_p} H^0(V_{K_p},\Of_{\X^*_{K_p K^p}})\fl H^0(V,\Of_{\X^*_{K^p}}).\]
\end{enumerate}
\end{enumerate}

\end{theo}

De plus, le th\'eor\`eme ci-dessus s'\'etend \`a toutes les vari\'et\'es
de Shimura de type Hodge (th\'eor\`eme IV.1.1 de \cite{S}), en particulier
aux vari\'et\'es de Shimura de type PEL.

La pr\'esence du bord de $\X^*_K$ cause quelques probl\`emes techniques (dont
la plupart sont trait\'es dans les sections II.2 et II.3 de \cite{S}).
Nous allons
donner une esquisse de preuve qui ignore ces probl\`emes.

On fixe $K^p\subset\G(\Af^p)$, et on note, pour tout $m\in\Nat^*$,
\[\Gamma_0(p^m)=\{\gamma\in\GSp_{2n}(\Z_p)|\gamma=\begin{pmatrix}* & * \\
0 & * \end{pmatrix}\mod p^m\mbox{ et }\det(\gamma)=1\mod p^m\},\]
\[\Gamma(p^m)=\{\gamma\in\GSp_{2n}(\Z_p)|\gamma=\begin{pmatrix}1 & 0 \\
0 & 1 \end{pmatrix}\mod p^m\}.\]

Dans la premi\`ere partie de la preuve, on prouve (i) pour la limite sur
les $K_p=\Gamma_0(p^m)$ et dans un voisinage du lieu anticanonique. Expliquons
ce que ceci signifie.

Soit $A\fl S$ un sch\'ema ab\'elien sur un sch\'ema $S$ de caract\'eristique $p$. Si $A^{(p)}$ est le changement de base de $A\fl S$ par le morphisme de
Frobenius absolu de $S$, alors il existe une unique isog\'enie $V:A^{(p)}\fl
A$ (appel\'ee \emph{Verschiebung}) dont le compos\'e (dans les deux sens)
avec le morphisme de
Frobenius relatif $A\fl A^{(p)}$ est la multiplication par $p$.
Cette isog\'enie induit un morphisme $V^*:\omega_{A/S}\fl\omega_{A^{(p)}/S}=
\omega_{A/S}^{\otimes p}$,\footnote{$\omega_{A/S}=\wedge^{\dim(A/S)}(\Omega^1_{A/S})$.}
c'est-\`a-dire une section globale $\Ha(A/S)$
de $\omega_{A/S}^{\otimes(p-1)}$, appel\'ee \emph{invariant de Hasse de $A$}.
Rappelons que $\Ha(A/S)$ est inversible si et seulement si $A$ est ordinaire
\footnote{C'est-\`a-dire que $A[p](x)$ a exactement $p^{\dim(A/S)}$ \'el\'ements
pour tout point g\'eom\'etrique $x$ de $S$.} (lemme III.2.5 de \cite{S}).
L'invariant de Hasse permet donc de mesurer \`a quel point $A$ est loin
d'\^etre ordinaire.

Dans la suite, si on \'ecrit $p^\varepsilon$,
avec $\varepsilon\in[0,1[$, on supposera
toujours que $\varepsilon$ est dans l'image de $\Z_p^{\rm cycl}$ par la valuation
$p$-adique, et $p^\varepsilon$ sera un \'el\'ement (quelconque) de $\Z_p^{\rm cycl}$
de valuation $\varepsilon$.

On note $\X=\X_{\G(\Z_p)K^p}$ (la vari\'et\'e de Shimura de niveau trivial en $p$).
Pour $\varepsilon$ comme ci-dessus, on d\'efinit un ouvert $\X(\varepsilon)$
par la condition $|\Ha|\geq |p^\varepsilon|$ (o\`u on prend l'invariant de Hasse
de la vari\'et\'e ab\'elienne universelle). Pour que cela ait un sens, on doit
d'abord d\'efinir $\X(\varepsilon)$ sur un mod\`ele entier de $\X$ (qui existe
car on a pris le niveau trivial en $p$) en utilisant un rel\`evement local
de l'invariant de Hasse, et montrer que le r\'esultat ne d\'epend pas de
ce rel\`evement; voir la d\'efinition III.2.11 et le lemme III.2.12
de \cite{S}. On note aussi $\Abf\fl\X$ l'espace adique associ\'e au
sch\'ema ab\'elien universel,
et $\Abf(\varepsilon)=\Abf\times_\X\X(\varepsilon)$. Si $K_p\subset\G(\Q_p)$ est
quelconque, on note $\X_{K_p K^p}(\varepsilon)=\X_{K_p K^p}\times_\X\X(\varepsilon)$.

On fixe $\varepsilon$ et $m\in\Nat$. Alors, sur les mod\`eles entiers, le
sch\'ema ab\'elien universel $\Abf(p^{-m}\varepsilon)\fl\X(p^{-m}\varepsilon)$
admet un sous-groupe canonique de niveau $m$, c'est-\`a-dire un sous-groupe
$C_m$ du groupe des points de $p^m$-torsion qui est \'egal \`a $\Ker(F^m)$
modulo $p^{1-\varepsilon}$ ($F$ est le Frobenius absolu); voir le (ii) du
th\'eor\`eme III.2.14 de \cite{S}. Comme $X_{\Gamma_0(p^m)\K^p}$ param\`etre les
vari\'et\'es ab\'eliennes (principalement polaris\'ees) $A$
munies d'une structure de niveau $K^p$ et d'un sous-groupe totalement 
isotrope de
$A[p^m]$, l'existence de $C_m$ donne un morphisme $\X(p^{-m}\varepsilon)
\fl\X_{\Gamma_0(p^m)}$. Scholze prouve alors les r\'esultats suivants (th\'eor\`eme
III.2.14(iii),(iv) et lemme III.2.16 de \cite{S}) :
\begin{enumerate}
\item[(a)] Pour tout $m\geq 1$, on a un diagramme cart\'esien
\[\xymatrix{\X(p^{-m-1}\varepsilon)\ar[r]\ar[d] & \X_{\Gamma_0(p^{m+1})K^p}\ar[d] \\
\X(p^{-m}\varepsilon)\ar[r] & \X_{\Gamma_0(p^{m})K^p}}\]
o\`u les fl\`eches horizontales sont celles d\'efinies ci-dessus, la fl\`eche
verticale de droite est la projection canonique et la fl\`eche verticale de
gauche est un rel\`evement du Frobenius.

\item[(b)] L'image $X_{\Gamma_0(p^m)K^p}(\varepsilon)_a$ du morphisme
$\X(p^{-m}\varepsilon)\fl\X_{\Gamma_0(p^m)}$ est un sous-espace ouvert et ferm\'e
de $\X_{\Gamma^0(p^m)K^p}(\varepsilon)$, d\'efini par la condition que
$D[p]\cap C_1=\{0\}$, o\`u $D$ est le sous-groupe totalement isotrope donn\'e
par la structure de niveau $\Gamma_0(p^m)$.\footnote{Le {\og a\fg} signifie
{\og anticanonique\fg}.}
De plus,
$X_{\Gamma_0(p^m)K^p}(\varepsilon)_a$ est affino\"ide pour $m$ assez grand.
\end{enumerate}

On d\'eduit des r\'esultats ci-dessus qu'il existe un espace perfecto\"ide
affino\"ide $\X_{\Gamma_0(p^\infty)}(\varepsilon)_a$ tel que
\[\X_{\Gamma_0(p^\infty)}(\varepsilon)_a\sim\varprojlim_m\X_{\Gamma_0(p^m)K^p}
(\varepsilon)_a.\]

La deuxi\`eme \'etape consiste \`a passer au niveau $\Gamma(p^m)$, afin
d'obtenir un espace affino\"ide perfecto\"ide $\X_{\Gamma(p^\infty)}
(\varepsilon)_a$ tel que
\[\X_{\Gamma(p^\infty)}(\varepsilon)_a\sim\varprojlim_m\X_{\Gamma(p^m)K^p}
(\varepsilon)_a.\]
La seule difficult\'e vient du fait que les morphismes $\X^*_{\Gamma(p^m)K^p}\fl
\X^*_{\Gamma_0(p^m)K^p}$ ne sont pas \'etales au bord. Voir la section III.2.5
de \cite{S} pour les d\'etails.

Pour tout espace adique $\Y$, on note $|\Y|$ l'espace topologique sous-jacent.
Soit
\[|X_{\Gamma(p^\infty)}|=\varprojlim_m|\X_{\Gamma(p^m)K^p}|.\]
On a trouv\'e un ouvert $|\X_{\Gamma(p^\infty)}(\varepsilon)_a|$
de cet espace dont chaque point admet un voisinage
(venant d'un) perfecto\"ide affino\"ide. Notons que le groupe $\G(\Q_p)$ agit
de mani\`ere continue
sur $|\X_{\Gamma(p^\infty)}|$, et que la condition {\og avoir un voisinage
perfecto\"ide affino\"ide\fg} est stable par cette action. L'id\'ee est maintenant
de montrer que les $\G(\Q_p)$-translat\'es de $|\X_{\Gamma(p^\infty)}(\varepsilon)_a|$ recouvrent $|\X_{\Gamma(p^\infty)}|$. Pour cela, Scholze utilise le morphisme
des p\'eriodes de Hodge-Tate.

Si $A$ est une vari\'et\'e ab\'elienne sur une extension compl\`ete et
alg\'ebriquement close $C$ de $\Q_p$, la suite spectrale de Hodge-Tate (voir par exemple
le th\'eor\`eme 3.20 de \cite{S2}) donne une suite exacte courte
\[0\fl (\Lie A)(1)\fl T_p A\otimes_{\Z_p} C\fl (\Lie A^*)^*\fl 0,\]
o\`u $T_p A$ est le module de Tate $p$-adique de $A$. Si $A$ vient d'un point
de $|\X_{\Gamma(p^\infty)}|$, alors on a un isomorphisme $T_p A\simeq\Z_p^{2n}$
(donn\'e par la structure de niveau infini en $p$), donc la suite exacte
ci-dessus d\'efinit un point de $\Fl(C)$. On obtient ainsi une application
$\G(\Q_p)$-\'equivariante $|\pi_{HT}|:|\X_{\Gamma(p^\infty)}|\fl|\Fl|$, dont on
montre qu'elle est continue en regardant ce qui se passe sur des voisinages
pro-\'etales perfecto\"ides des points (lemme III.3.4 de \cite{S}). En effet,
la trivialisation du module de Tate $p$-adique $T_p A$ donn\'ee par la
structure de niveau infinie en $p$ existe en fait sur des voisinages
pro-\'etales des points de la vari\'et\'e de Shimura perfecto\"ide,
\footnote{Il est ici important de disposer de la topologie pro-\'etale. En
effet, sur un voisinage \'etale, on ne pourrait obtenir que des
trivialisations des groupes de points de torsion $A[p^N]$.}
ce qui permet, sur un tel voisinage,
de d\'efinir une filtration de Hodge-Tate relative, et donc de montrer que
le morphisme de Hodge-Tate est un morphisme d'espaces adiques (et en particulier
continu).

On a le plongement de Pl\"ucker $\Fl\subset\Proj^{\left(\substack{2n \\ n}\right)-1}$
(d\'efini en envoyant $W$ sur $\wedge^n W$). On note $s_J$ les coordonn\'ees
homog\`enes sur le but de ce plongement, o\`u $J$ parcourt les sous-ensembles
de $\{1,\dots,2n\}$ de cardinal $n$. Si on fixe
un tel $J$, on a un ouvert affino\"ide
$\Fl_J$ de $\Fl$, d\'efini par la condition $|s_{J'}|\leq |s_J|$, pour tout
$J'$.

En utilisant les r\'esultats de son article \cite{SW} avec Weinstein (plus
pr\'ecis\'ement, le th\'eor\`eme B de cet article), Scholze montre que
l'image inverse de $\Fl(\Q_p)$ par $|\pi_{HT}|$ est le lieu ordinaire
(lemme III.3.6 de \cite{S}), puis qu'il existe un voisinage ouvert $U$ du
point $x=0^n\oplus\Q_p^n$ de $\Fl$ tel que $|\pi_{HT}|^{-1}(U)\subset
|\X_{\Gamma_0(p^\infty)}(\varepsilon)_a|$ (pour $\varepsilon>0$ convenable).
Or, si $\gamma\in\G(\Q_p)$
est l'\'el\'ement diagonal $(p,\dots,p,1,\dots,1)$, alors $\gamma^N
\Fl_{\{n+1,\dots,2n\}}$ pour $N$ assez grand. Comme $\G(\Z_p).\Fl_{\{n+1,\dots,2n\}}=
\Fl$, on en d\'eduit que $\G(\Q_p).U=\Fl$, ce qui permet de construire
un {\og atlas perfecto\"ide\fg} de $|\X_{\Gamma(p^\infty)}|$, donc d'obtenir l'espace
perfecto\"ide $\X_{\Gamma(p^\infty)}$  (corollaire III.3.11 de \cite{S}). Une
fois que l'on a $\X_{\Gamma(p^\infty)}$, il est assez facile de montrer que
$|\pi_{HT}|$ vient d'un morphisme d'espaces adiques
$\pi_{HT}:\X_{\Gamma(p^\infty)}\fl\Fl$ (et un peu moins facile de montrer que ce
morphisme s'\'etend au bord, voir les corollaires III.3.12 et III.3.16 de
\cite{S}).

\section{Cohomologie compl\'et\'ee et faux invariants de Hasse}

\subsection{Cohomologie compl\'et\'ee}\label{cc}

On utilise les notations de la section pr\'ec\'edente, et on fixe
toujours un sous-groupe ouvert compact (assez petit) $K^p$ de $\G(\Af^p)$.
Rappelons que la cohomologie compl\'et\'ee de la vari\'et\'e de Shimura de
$\G$ en niveau mod\'er\'e, dont la d\'efinition est due
\`a Calegari et Emerton (\cite{CE}), est donn\'ee par la formule :
\[\widetilde{H}^i_{c,K^p}=\varprojlim_m\varinjlim_{K_p} H^i_c(X^{\G}_{K_p K^p},
\Z/p^m\Z).\]

C'est un $\Z_p$-module $p$-adiquement complet, qui admet une action
de $\G(\Q_p)\times\Hf^S$, pour tout ensemble fini de nombres premiers $S$
tel que $K^p$ soit hypersp\'ecial en dehors de $S$.\footnote{Bien s\^ur,
$H^i_{c,K^p}$ admet en fait une action de toute l'alg\`ebre de Hecke mod\'er\'ee
de niveau $K^p$.} On fixe un tel $S$.

\begin{rema} La cohomologie compl\'et\'ee permet de d\'efinir
une notion assez g\'en\'erale
de {\og repr\'esentation automorphe $p$-adique\fg} pour un groupe r\'eductif
quelconque (c'\'etait d'ailleurs l'une des motivations de Calegari et Emerton).

\end{rema}

On note aussi, pour tout $m\in\Nat$,
\[\widetilde{H}^i_{c,K^p}(\Z/p^m\Z)=\varinjlim_{K_p} H^i_c(X^{\G}_{K_p K^p},
\Z/p^m\Z).\]

Les r\'esultats des sections \ref{comp} et \ref{HT} permettent de montrer sans
trop de peine le r\'esultat suivant :

\begin{theo}[Th\'eor\`eme IV.2.1 de \cite{S}]\label{th:C} Soit $\Ical\subset
\Of_{\X^*_{\Gamma(p^\infty)K^p}}$ l'id\'eal du bord, et $\Ical^+=\Ical\cap
\Of_{\X^*_{\Gamma(p^\infty)K^p}}^+$. Soit $C$ une extension compl\`ete et
alg\'ebriquement close de $\Q_p$.

Alors on a des presque isomorphismes naturels
\[\widetilde{H}^i_{c,K^p}(\Z/p^m\Z)\otimes_{\Z/p^m\Z}\Of_C^a/p^m\simeq
H^i(\X^*_{\Gamma(p^\infty)K^p},\Ical^{+a}/p^m),\]
compatibles avec l'action de $\Hf^S$.\footnote{Voir l'\'enonc\'e du th\'eor\`eme
IV.2.1 de \cite{S} pour la d\'efinition de cette action.}.

\end{theo}

\begin{rema} En particulier, on peut voir la cohomologie compl\'et\'ee comme
la cohomologie \'etale de la vari\'et\'e de Shimura perfecto\"ide, ce qui
en donne une interpr\'etation naturelle.

\end{rema}

\subsection{Faux invariants de Hasse}\label{faux}

On a presque fini la preuve du th\'eor\`eme \ref{th:princ}, ou de sa version
plus pr\'ecise, le th\'eor\`eme IV.3.1 de \cite{S}, qui dit, en gros,
que, si on fixe $m\geq 1$, alors
tout caract\`ere de $\Hf^S$ qui appara\^it dans un $\widetilde{H}^i_{c,
\K^p}(\Z/p^m\Z)$ appara\^it aussi dans un $H^0(\X^*_{K_p K^p},\omega^{mk}\otimes
\Ical)$, pour un $K_p\subset\G(\Q_p)$ compact ouvert (variable) et
un entier $k\geq 1$ (variable aussi). Ici, $\Ical$ d\'esigne comme en niveau
infini l'id\'eal du bord.

Pour finir la preuve, il faut pouvoir passer des groupes 
$H^i(\X^*_{\Gamma(p^\infty)K^p},\Ical^+/p^m)$ (qui sont en niveau infini et en
degr\'e cohomologique quelconque) aux groupes
\mbox{$H^0(\X^*_{K_p K^p},\omega^{mk}\otimes\Ical)$.}

On utilise le recouvrement de $\X^*_{\Gamma(p^\infty)K^p}$ par les ouverts
$\Vcal_J:=(\pi_{HT})^{-1}(\Fl_J)$ donn\'e par le (iv)
du th\'eor\`eme \ref{th:princ}. Comme les $\Vcal_J$ sont
affino\"ides perfecto\"ides (et gr\^ace \`a une propri\'et\'e technique
du bord, voir le (ii) du th\'eor\`eme IV.1.1 de \cite{S}), la cohomologie
de $\Ical^+/p^m$ sur ces ouverts est presque concentr\'ee en degr\'e $0$.
De plus, comme $\pi_{HT}$ est \'equivariant pour les op\'erateurs de Hecke
en dehors de $p$ (agissant trivialement sur $\Fl$), les ouverts $\Vcal_J$
sont stables par ces op\'erateurs, donc $\Hf^S$ agit encore sur les
$H^i(\Vcal_J,\Ical^+/p^m)$.

On utilise la deuxi\`eme partie du point (iv) du th\'eor\`eme \ref{th:princ}
pour montrer que tout caract\`ere de $\Hf^S$ apparaissant dans un
$\H^0(\Vcal_J,\Ical^+/p^m)$ appara\^it en fait dans un
$\H^0(\Vcal_{J,K_p},\Ical^+/p^m)$, o\`u $K_p\subset\G(\Q_p)$ est assez petit et
$\Vcal_{J,K_p}\subset\X^*_{K_p K^p}$ est un ouvert affino\"ide d'image inverse
$\Vcal_J$ dans $\X^*_{\Gamma(p^\infty)K^p}$.

Il faut encore montrer comment \'etendre les sections de $\Ical^+/p^m$ sur
$\Vcal_{J,K_p}$ qui sont vecteurs propres pour $\Hf^S$
\`a $\X^*_{K_p K^p}$ tout entier sans changer les valeurs
propres. La m\'ethode classique consiste \`a multiplier ces sections
propres par une puissance assez grande de l'invariant de Hasse. Ici,
on utilise plut\^ot les {\og faux invariants de Hasse\fg}, qui sont des
\'el\'ements de $H^0(\Vcal_{J,K_p},\omega)$ obtenus par pullback de sections
bien choisies dans $H^0(\Fl_J,\omega_\Fl)$ (et par descente \`a un niveau fini
assez petit $K_p$); voir le lemme II.1.1 et la page 72 de \cite{S}.
Le fait que
la multiplication par ces faux invariants de Hasse ne
change pas les valeurs propres de $\Hf_S$ r\'esulte de la propri\'et\'e
d'\'equivariance de~$\pi_{HT}$.

\section{Quelques autres applications}\label{appl}

Indiquons deux autres applications des r\'esultats de \cite{S}. (Cette liste
d'applications ne se veut en aucun cas exhaustive.)

\subsection{Cohomologie compl\'et\'ee}

Le th\'eor\`eme \ref{th:C} donne une formule pour la cohomologie compl\'et\'ee
de la section \ref{cc}. En utilisant ce th\'eor\`eme et le fait que les
espaces topologiques sous-jacents aux $\X^*_{K_p K^p}$ sont de dimension
cohomologique $\leq d$, o\`u $d=n(n+1)/2$ est la dimension des vari\'et\'es
alg\'ebriques $X_{K_p K^p}$, Scholze en d\'eduit que $\widetilde{H}^i_{c,K^p}=0$
pour $i>d$, puis une grande partie de la conjecture 1.5 de \cite{CE} (corollaire
IV.2.3 de \cite{S}).
\footnote{L\`a encore, les r\'esultats cit\'es sont en fait vrais pour toutes
les vari\'et\'es de Shimura de type Hodge, en particulier celles de type PEL.}

\subsection{Mod\`eles entiers \'etranges}

Les faux invariants de Hasse de la section \ref{faux} permettent de d\'efinir
des mod\`eles entiers jusqu'ici inconnus des vari\'et\'es de Shimura de type
Hodge. Dans le preprint \cite{PSt}, Pilloni et Stroh ont \'etudi\'e ces mod\`eles
entiers et les ont utilis\'es pour construire des repr\'esentations galoisiennes
associ\'ees \`a des repr\'esentations automorphes non n\'ecessairement
cohomologiques\footnote{Mais apparaissant dans la cohomologie coh\'erente
de fibr\'es vectoriels automorphes.}
du groupe $\G$ d\'efinissant les vari\'et\'es de Shimura.

\end{document}